\documentclass[a4paper,11pt]{article}
\usepackage{amsmath,amsthm,amssymb}

\renewcommand{\baselinestretch}{1}

\makeatletter
  
   \@addtoreset{equation}{section}
\makeatother
\def\<{\left<} \def\>{\right>}
\newtheorem{theorem}{Theorem}
\newtheorem{lemma}[theorem]{Lemma}

\def\bea{\begin{eqnarray} }
\def\eea{\end{eqnarray} }
\def\be{\begin{equation} }
\def\ee{\end{equation} }
\def\qed{\ifhmode\unskip\nobreak\fi\ifmmode\ifinner\else\hskip5pt
\fi\fi\hbox{\hskip5 pt \vrule width4 pt height6 pt depth1.5 pt \hskip1pt }}

\begin{document}
\title{Legendre surfaces with harmonic mean curvature vector field in the unit $5$-sphere
\footnote{Rocky Mountain Journal of Mathematics  {\bf 40} (2010), 313-320. An appendix 
is added.}}
\author{Toru Sasahara}

\date{}
\maketitle

\begin{abstract}
{\footnotesize  
 We obtain the explicit representation of Legendre surfaces in the unit $5$-sphere
with harmonic mean curvature vector field, under the condition that the mean curvature
function is constant along a certain special direction.}
\end{abstract}


\section{Introduction}
It is well-known that an odd-dimensional unit sphere $S^{2n+1}$ is equipped with the standard Sasakian structure 
$(g, \phi, \eta, \xi)$ (see, Section 2). The study of minimal Legendre submanifolds in $S^{2n+1}$ is a very 
active field and closely related  to one of minimal Lagrangian submanifolds in complex projective space.

A natural generalization of a minimal submanifold is 
a submanifold with parallel mean curvature vector field.
However, Legendre submanifolds in $S^{2n+1}$ with parallel mean curvature vector
 field are in fact minimal (see, \cite{yk}).
So, as a generalization of minimal Legendre submanifolds,  it is natural to consider 
Legendre submanifolds whose mean curvature vector field $H$ is harmonic with respect
to the normal Laplacian $\Delta^D$ (see, Section 2), that is, 
\bea
\Delta^D H=0.
\eea

The purpose of this paper is to study the class of 
nonminimal Legendre surfaces satisfying $(1.1)$ in the unit $5$-sphere.
The mean curvature function $||H||$ of such surfaces is not constant (see, Lemma 2 in Section 3). 
From this, the class seems to be very large. Thus, we need to add some natural condition. 
On nonminimal Legendre submanifolds,
there exists a special vector field: $\phi H$.
Moreover in case the dimension is $2$, up to signs, there is a unique unit vector field 
normal to $\phi H$. We denote it by $(\phi H)^{\perp}$.
In this paper, under the condition that the square of 
the mean curvature function is constant along one of $\phi H$ and $(\phi H)^{\perp}$,
we completely determine
nonminimal Legendre surfaces satisfying $(1.1)$ in the unit $5$-sphere as follows:

\begin{theorem}
Let $f: M^2\rightarrow S^5\subset{\bf C}^3$ be a nonminimal Legendre 
 surface satisfying $\Delta^D H=0$
in the unit $5$-sphere. Then the mean curvature function of $M^2$ is not constant.

\noindent If $H$ satisfies $\phi H||H||^2=0$, then $f$ is locally given by
\be
f(x,y)=\left(\frac{1}{\sqrt{2}}{\rm exp}(\frac{1+\sqrt{5}}{2}iy){
\rm cos}x, \frac{1}{\sqrt{2}}{\rm exp}(\frac{1-\sqrt{5}}{2}iy){
\rm cos}x, {\rm sin}x\right).
\ee

\noindent If $H$ satisfies $(\phi H)^{\perp}||H||^2=0$,  then
$f$ is locally given by
\bea
f(x,y)=\frac{1}{\sqrt{2}}\Biggl(i+{\rm sin}x, ({\rm sec}x+{\rm
tan}x)^i{\rm cos}x{\rm cos}y, ({\rm sec}x+{\rm
tan}x)^i{\rm cos}x{\rm sin}y\Biggr).
\eea

\end{theorem}

\section{Legendre submanifolds in the unit sphere}
Let ${\bf C}^{n+1}$ be the complex Euclidean $(n+1)$-space 
together with the canonical complex structure $J$. 
Denote by $S^{2n+1}$ the unit sphere with the standard induced metric $g$ in ${\bf C}^{n+1}$. 
The position vector field $\mathbf{x}$ of 
$S^{2n+1}$ is a unit normal vector field of
$S^{2n+1}$ in ${\bf C}^{n+1}$
and  the vector field $\xi:=-J\mathbf{x}$ is tangent to $S^{2n+1}$.
Define a $1$-form $\eta$ and an endomorphism field $\varphi$ on 
$M$ by the formula:
$$
JX=\phi X+\eta(X){\bf x}, X \in TM.
$$ 
It is easy to see that $(g, \varphi, \eta, \xi)$ satisfies
$$
\varphi^{2}=-I+\eta \otimes \xi,\ d\eta(X,Y) =g(X,\varphi Y).
$$
Thus $(g, \varphi,\eta,\xi)$ is a contact metric structure of $S^{2n+1}$.

Denote by $\overline{\nabla}$ the Levi-Civita connection of 
$S^{2n+1}$. Then 
$$
{\overline \nabla}_{X}\xi=-\varphi X,\ 
(\overline{\nabla}_{X}\varphi)Y=g(X,Y)\xi-\eta(Y)X.
$$
These formulas imply that $S^{2n+1}$ is a Sasakian  manifold.

An immersed $n$-submanifold $x:M^{n}\to S^{2n+1}$ is said to be a {\it Legendre submanifold} if
$x^{*}\eta=0$.
 The
formulas of Gauss and Weingarten of $x$ are given respectively by
\bea
\begin{split}
&\bar \nabla_XY= \nabla_XY+h(X,Y),\\
&\bar\nabla_X V = -A_{V}X+D_XV,
\end{split}
\eea
where $X, Y\in  TM^n$, $V\in T^{\perp}M^n$, 
$\nabla$, $h,A$ and $D$ are the Levi-Civita connection of $M^n$, the second fundamental
form, the shape operator and the normal
connection. The mean curvature vector $H$ is given by $H=\frac{1}{n}{\rm trace}h$.
Its length $||H||$ is called the {\it mean curvature function} of $M^n$.
The  normal Laplacian 
is defined by $\Delta^D=-\sum_{i=1}^{n}
(D_{e_i}D_{e_i}-D_{\nabla_{e_i}e_i})$, where
$\{e_i\}$ is a local orthonormal frame  of $M^n$.

For Legendre submanifolds we have \cite{bl}
\be
A_{\phi Y}X=-\phi h(X, Y)=A_{\phi X}Y,\quad A_{\xi}=0.
\ee
Moreover a straightforward computation shows that 
the equations of Gauss, Codazzi, Ricci of Legendre submanifolds in the unit sphere
are equivalent to
\bea
& &\<R(X,Y)Z,W\> =  \<[A_{\phi Z}, A_{\phi W}]X,Y\>+\<X, W\>\<Y, Z\>-\<X, Z\>\<Y, W\>,\\
& & ({\bar\nabla}_{X}h)(Y,Z)=({\bar\nabla}_{Y}h)(X,Z),
\eea
where $\bar\nabla h$ is defined by
 $({\bar\nabla}_{X}h)(Y,Z)= D_X h(Y,Z) - h(\nabla_XY,Z) - h(Y,\nabla_X Z)$.

\section{The proof of Theorem 1}
We assume that the mean curvature function is nowhere zero.
Let $\{e_i\}$ ($i=1,\ldots, 5$) be an orthonormal frame along $M^2$ such that 
$e_1, e_2$ are tangent to
$M^2$, $\phi e_1=e_3$, $\phi e_2=e_4$, $\xi=e_5$ and $H=\frac{\alpha}{2}\phi e_1$, with $\alpha>0$. 
Then, it follows from (2.1) and (2.2) that the
 second fundamental form takes the form:
\bea
&&h(e_1, e_1)=(\alpha-c)\phi e_1+b\phi e_2,\nonumber\\
&&h(e_1, e_2)=b\phi e_1+c\phi e_2,\\ 
&&h(e_2, e_2)=c\phi e_1-b\phi e_2,\nonumber
\eea
 for some functions $b$, $c$.  

We put $\omega_i^j(e_k)=\<\nabla_{e_k}e_i, e_j\>$. From (2.4) we get 
\bea
& &e_1c+3b\omega_1^2(e_1)=e_2b+(\alpha-3c)\omega^2_1(e_2),\\
& &-e_1b+3c\omega_1^2(e_1)=e_2c+3b\omega^2_1(e_2),\\
& &e_2(\alpha-c)-3b\omega_1^2(e_2)=e_1b+(\alpha-3c)\omega^2_1(e_1).
\eea

Suppose that $M^2$ satisfies $\Delta^DH=0$.
Then we have the following  three lemmas (see p.290 and 291 in \cite{sa}).
\begin{lemma}
\bea
& &\Delta_M\alpha+\alpha\{1+(\omega^2_1(e_1))^2
+(\omega^2_1(e_2))^2\}=0,\\
& &2(e_1\alpha)\omega^2_1(e_1)+2(e_2\alpha)\omega^2_1(e_2)
+\alpha\{e_1(\omega^2_1(e_1))+e_2(\omega^2_1(e_2))\}=0,\\
& &e_1\alpha+\alpha\omega^2_1(e_2)=0,
\eea
where $\Delta_M$ is the Laplace operator acting on $C^{\infty}(M)$.
\end{lemma}
Note that $\alpha$ is not constant by (3.5), since $\alpha$ is nowhere zero.

 \begin{lemma}
There exists local coordinates $x$, $y$ such that
\bea 
&&e_1=\alpha\partial_x, \nonumber\\
&&e_2=\alpha\partial_y,\nonumber\\
&&\omega_1^2(e_1)=\alpha_y,\nonumber\\ 
&&\omega_1^2(e_2)=-\alpha_x.\nonumber
\eea
\end{lemma}
\begin{lemma}The following relation holds:
\bea 
b^2=\frac{\alpha c}{2}-c^2.
\eea
\end{lemma}


The allied mean curvature vector  $a(H)$ is defined by 
$\sum_{r=4}^{5}({\rm trace}A_HA_{e_r})e_r$ (see, p.197 in \cite{ch1}). 
If $a(H)$ vanishes identically on $M^2$, 
it is called a  $Chen$ 
surface. 

Suppose that $M^2$ is not Chen surface, i.e., $b\ne0$. 
Then we may assume that $b>0$, if necessary by changing the sign of $e_2$.
By differentiating (3.8) we get
\be
b_i=\frac{(\alpha-4c)c_i+\alpha_i c}{4b},
\ee 
where $i=x, y$.

Using (3.9) we replace (3.2) and (3.3) by the derivatives with respect to $x$ and $y$ as follows:
\be
\left(\begin{array}{cccccc}
\alpha& -\dfrac{\alpha(\alpha-4c)}{4b}\\
\dfrac{\alpha(\alpha-4c)}{4b} & \alpha
\end{array}
\right)
\left(\begin{array}{cccccc}
c_x\\
c_y
\end{array}
\right)
=\left(\begin{array}{cccccc}
-(\alpha-3c)\alpha_x-\dfrac{12b^2-\alpha c}{4b}\alpha_y\\
\dfrac{12b^2-\alpha c}{4b}\alpha_x+3c\alpha_y
\end{array}
\right).
\ee

First we investigate the case of $\phi H||H||=0$, i.e., $\alpha_x=0$. 
Then using  (3.8) and (3.10) we get
\bea
&&c_x=-\frac{8bc}{\alpha^2}\alpha_y,\\
&&c_y=\frac{5\alpha c-8c^2}{\alpha^2}\alpha_y.
\eea
By a long but straightforward computation, we obtain
\bea 
&&bc_{xy}=
(-24\alpha_y^2-4\alpha\alpha_{yy})\biggl(\frac{c}{\alpha}\biggr)^2
+(112\alpha_y^2+8\alpha\alpha_{yy})\biggl(\frac{c}{\alpha}\biggr)^3\nonumber\\
&&\quad\quad-128\alpha_y^2\biggl(\frac{c}{\alpha}\biggr)^4,\\
&&bc_{yx}=
-20\alpha_y^2\biggl(\frac{c}{\alpha}\biggr)^2
+104\alpha_y^2\biggl(\frac{c}{\alpha}\biggr)^3
-128\alpha_y^2\biggl(\frac{c}{\alpha}\biggr)^4,
\eea
Since $bc_{xy}=bc_{yx}$, we find that
$\alpha_y^2+\alpha\alpha_{yy}=0$ or $c=0$ or $\alpha=2c$.
If $\alpha_y^2+\alpha\alpha_{yy}=0$, it contradicts to (3.5). If  $c=0$ or $\alpha=2c$,
it is a contradiction to $b\ne 0$ by (3.8). 

Next, we investigate the case of $(\phi H)^{\perp}||H||=0$, i.e., $\alpha_y=0$. 
Similarly as in  the case of  $\alpha_x=0$,  we have
\bea
&&c_x=\frac{-3\alpha c+8c^2}{\alpha^2}\alpha_x,\\
&&c_y=\frac{4\alpha b -8bc}{\alpha^2}\alpha_x,\\
&&bc_{xy}=-6\alpha_{x}^2\biggl(\frac{c}{\alpha}\biggr)
+56\alpha_x^2\biggl(\frac{c}{\alpha}\biggr)^2 \nonumber\\
&&\quad\quad-152\alpha_x^2\biggl(\frac{c}{\alpha}\biggr)^3
+128\alpha_x^2\biggl(\frac{c}{\alpha}\biggr)^4,\\
&&bc_{yx}=(-4\alpha_x^2+2\alpha\alpha_{xx})\biggl(\frac{c}{\alpha}\biggr)
+(48\alpha_x^2-8\alpha\alpha_{xx})\biggl(\frac{c}{\alpha}\biggr)^2\nonumber\\
&&\quad\quad-(144\alpha_x^2-8\alpha\alpha_{xx})\biggl(\frac{c}{\alpha}\biggr)^3
+128\alpha_x^2\biggl(\frac{c}{\alpha}\biggr)^4,
\eea
Hence we get $\alpha^2_x+\alpha\alpha_{xx}=0$ or $c=0$ or $\alpha=2c$,
since $bc_{xy}=bc_{yx}$. 
It is a contradiction as mentioned above.

Therefore we conclude that $b$ must be $0$, i.e., $M^2$ must be a  Chen surface if $\phi H||H||^2=0$
 or $(\phi H)^{\perp}|H||^2=0$. Applying the classification
of Legendre Chen surfaces satisfying $\Delta^DH=0$ (see, Theorem 8 and Corollary 9 in \cite{sa2}), 
we can prove the statement.

\section{Other examples of Legendre surfaces satisfying $\Delta^DH=0$ }

In this section, we show a way to construct Legendre surfaces satisfying 
$\Delta^DH=0$ and $(\phi H)^{\perp}||H||^2\ne0$ and $\phi H||H||^2\ne 0$.

One can obtain the following existence and uniqueness  theorem by the
similar way to those given in \cite{cv1} and \cite{cv2} (cf. \cite{cdvv}).
\begin{theorem}
Let $(M^n , \<\cdot,\cdot\>)$ be an n-dimensional simply connected Riemannian
manifold. Let $\sigma$ be a symmetric bilinear $TM^n$-valued
form on $M^n$ satisfying 

$(1)$ $\<\sigma(X, Y) ,Z\>$ is totally symmetric,

$(2)$ $(\nabla\sigma)(X, Y, Z)=\nabla_X\sigma(Y, Z)-\sigma(\nabla_XY, Z)
-\alpha(Y, \nabla_XZ)$ is totally symmetric,

$(3)$ $R(X, Y)Z=\<Y, Z\>X-\<X, Z\>Y
+\sigma(\sigma(Y, Z), X)-\sigma(\sigma(X,Z),Y)$.

\noindent
Then there exists  a Legendre isometric
immersion $x:(M^n , \<\cdot,\cdot\>)\rightarrow
S^{2n+1}$ such that the second fundamental form
$h$ satisfies $h(X, Y)=\phi\sigma(X,Y)$.
\end{theorem}
\begin{theorem}
Let $x^1, x^2: M^n\rightarrow S^{2n+1}$ be two Legendre
isometric immersions of a connected Riemannian $n$-manifold
into the unit sphere $S^{2n+1}$  with second fundamental forms
$h^1$ and $h^2$. If
$$\<h^1(X, Y), \phi x^1_{*}Z\>=\<h^2(X, Y), \phi x^2_{*}Z\>$$
for all vector fields $X$, $Y$, $Z$ tangent to $M^n$,
there exists an isometry $A$ of $S^{2n+1}$  such that
$x^1=A\circ x^2$.
\end{theorem}

Let $f(t)$ be a solution of the following ODE:
\be
\frac{d^2f}{dt^2}=\frac{1}{2}e^{-2f}.
\ee
We put $\alpha(x, y):=e^{f(x-y)}$. Let $(M^2, g=\frac{1}{\alpha^2}(dx^2+dy^2))$ be a Riemannian $2$-manifold.
We define a symmetric bilinear form $\sigma$ on $M^2$ by
\bea
&&\sigma(e_1, e_1)=\frac{3}{4}\alpha e_1+\frac{\alpha}{4}e_2,\nonumber\\
&&\sigma(e_1, e_2)=\frac{\alpha}{4}e_1+\frac{\alpha}{4}e_2,\\ 
&&\sigma(e_2, e_2)=\frac{\alpha}{4}e_1-\frac{\alpha}{4}e_2,\nonumber
\eea
where $e_1=\alpha\partial_x$ and $e_2=\alpha\partial_y$.
By a straight-forward computation, we find that $((M^2, g), \sigma)$ satisfies (1), (2) and (3)
of Theorem 5. Therefore there exists a unique Legendre surfaces in $S^5$ whose second fundamental form $h$
is given by $h=\phi\sigma$.
Moreover such a surface satisfies $\Delta^D H=0$ and 
$(\phi H)^{\perp}||H||^2\ne0$ and $\phi H||H||^2\ne 0$.

\vskip20pt
\noindent {\bf\Large
 Appendix:  Correction to Theorem 1} 
 {\normalsize(Added on December 10, 2014)}

\vskip10pt

\noindent As stated previously, Theorem 1 can be proved by applying  Theorem 8 and Corollary 9 in \cite{sa2}.
However, Corollary 9 in \cite{sa2} is incorrect. 
The immersion (44) in \cite{sa2} should be replaced by the immersion 
(4.9) in \cite{ch3}.
 Also, the immersion (45) in \cite{sa2}  should be replaced by
the immersion (4.6)  in 
\cite{ch2} such that $\lambda$ is constant.  

Accordingly, 
  (1.2) of this paper should be corrected as follows:
 \be
f(x,y)=\left(z(y){\cos}x, {\sin}x\right), \tag{a.1}
\ee
where  $z(y)$ is a unit speed Legendre curve with nonzero constant curvature 
in $S^3(1)$.

Also,  (1.3) should be corrected as follows:
\be
f(x,y)=(z_1(x), z_2(x)\cos y, z_2(x)\sin y),  \tag{a.2}
\ee
where  $(z_1, z_2)$ is a unit speed Legendre curve in $S^3(1)$ given by
\be
\begin{split}
&z_1=(ib+{\rm sin}x)/\sqrt{1+b^2},\\
&z_2=\exp(ib\ln({\sec}x+{\tan}x)){\cos}x/\sqrt{1+b^2} 
\end{split} \tag{a.3}
\ee
for some nonzero real number $b$. 

\vskip10pt
\noindent {\bf Remark A.1}
{\rm 
(i)  A  unit speed Legendre curve $z(y)$ of constant curvature $\lambda$ in $S^3(1)$ is congruent to
\bea
z(y)=(1/\sqrt{1+b^2})(e^{iby}, be^{-iy/b})\nonumber
\eea
with $b=(\lambda+\sqrt{\lambda^2+4})/2$.

\noindent (ii) The curvature of Legendre curve given by (a.3) in $S^3(1)$ 
is $b\sec x$.

\noindent (iii) Both surfaces given by (a.1) and (a.2)  have constant curvature one.}

{\renewcommand{\baselinestretch}{1}

\noindent Department of System and Information Engineering

\noindent Hachinohe Institute
of Technology

\noindent Hachinohe 031-8501, JAPAN

\noindent E-mail address:  {\tt sasahara@hi-tech.ac.jp}

\end{document}